\documentclass[5p, number, sort&compress, times, lefttitle]{elsarticle}

\journal{Elsevier}

\usepackage[utf8]{inputenc}
\usepackage{amsmath,amsfonts}
\usepackage{gensymb}
\usepackage{hyphenat}
\usepackage{subcaption}
\usepackage{multirow}
\usepackage{graphicx}
\usepackage{bm}
\usepackage{cancel}

\usepackage{tikz}
\usetikzlibrary{patterns,hobby}

\usepackage[colorlinks=true]{hyperref}

\usepackage{xcolor}

\newcommand{\fref}[1]{Fig.~\ref{#1}}

\newcommand{\vm}[1]{\bm{#1}}
\newcommand{\vx}{\vm{x}}
\newcommand{\U}[1]{\mathbb{U}_{#1}}

\renewcommand{\Re}{{\mathbb{R}}}

\graphicspath{ {./figs/} }

\graphicspath{ {./figs/} }

\begin{document}

\title{Virtual element method for modeling the deformation of multiphase composites} 

\author[1]{N.\ Sukumar\corref{cor1}}
\ead{nsukumar@ucdavis.edu}

\author[1]{John E. Bolander}
\ead{jebolander@ucdavis.edu}

\cortext[cor1]{Corresponding author}

\address[1]{Department of Civil and Environmental Engineering, 
University of California, Davis, CA 95616, USA}

\begin{abstract}
In this paper, we study applications of the 
virtual element method (VEM) 
for simulating the deformation of multiphase composites. 
The VEM is a Galerkin approach that is applicable to meshes that consist of arbitrarily-shaped polygonal and polyhedral (simple and nonsimple) elements. In the VEM,
the basis functions are defined as the solution of a local
elliptic partial differential equation, and are never explicitly
computed in the implementation of the method. 
The stiffness matrix of each element is built by using the elliptic projection operator of the internal virtual work (bilinear form) and it consists of two terms: a consistency term that is exactly computed (linear patch test is satisfied) and a correction term (ensures stability) that is orthogonal to affine displacement fields and has the right scaling. The VEM simplifies
mesh generation for a multiphase composite: a stiff inclusion  
can be modeled using a single polygonal or polyhedral element. Attributes of the virtual
element approach are highlighted through comparisons with Voronoi-cell lattice models, which provide discrete representations of material structure. The comparisons involve a suite of two-dimensional linear elastic 
problems: patch test, axisymmetric circular inclusion problem,  and the deformation of a three-phase composite. The simulations demonstrate the accuracy and flexibility of the virtual
element method.
\end{abstract}


\maketitle

\section{Introduction}\label{sec:intro}
Opportunities exist for designing multiphase
materials with improved composite properties~\cite{Ashby:2013}. In many cases, these materials consist of one or more dispersed particulate (or fibrous) phases within a binding phase. Along with the properties of the individual phases, typically the
behavior of the phase interfaces 
has primary influences on the composite properties, notably those related to fracture and mass transport. 
Concrete, which consists of aggregate inclusions embedded in a cement-based matrix, is a prime example of a multiphase particulate material that benefits from mesoscale analysis and design. 

Even though 
continuum approaches, including the finite element method, have been used for mesoscale modeling of concrete materials, various discrete modeling approaches have also received much interest. Particle-based lattice models are advantageous in the simple and natural way cracks and other forms of displacement discontinuity are 
represented~\cite{Cusatis:2011,Bolander:2021}, largely avoiding the stress-locking phenomenon associated with ordinary continuum representations of fracture.
Such lattice models permit deformation and fracture of inclusions (heterogeneities) to be efficiently represented and captured in simulations but they cannot in general exactly represent homogeneous deformation states (elastic homogeneity) for arbitrary Poisson's ratio $\nu$~\cite{Schlangen:1996}. In contrast, finite elements satisfy the patch test but the need for high-quality meshes for heterogeneous microstructures and the computational costs that are incurred limit the  number of inclusions that can be explicitly modeled. 

In this paper, we demonstrate the flexibility and capabilities that the
virtual element method (VEM)~\cite{Beirao:2013:BPV} affords to 
model the deformation of multiphase composites, such as cement-based materials that contain aggregate inclusions. Some of the previous  contributions in the modeling of the concrete mesostructure using the 
VEM are due to Benedetto et al.~\cite{Benedetto:2018:VEZ} and
Rivarola et al.~\cite{Rivarola:2019:MAV,Rivarola:2020:VEI}. We compare the performance of the VEM against a Voronoi-cell lattice model (VCLM) based on the rigid-body-spring concept of Kawai~\cite{Kawai:1978,Bolander:1998}. Our emphasis in this paper is to promote 
VEM as a methodology that has the desirable attributes of such lattice models as well as
the FEM to model multiphase materials (e.g., concrete composites).

\smallskip
The virtual element element (VEM)~\cite{Beirao:2013:BPV} is a stabilized 
high-order Galerkin discretizations on polygonal and polyhedral meshes to 
solve boundary-value problems. It provides a variational framework for the first order mimetic finite-difference scheme~\cite{Beirao:2014:MFD}, and is a generalization of hourglass finite elements~\cite{Flanagan:1981:AUS} to polytopal meshes~\cite{Cangiani:2015:HSV}.
In the VEM, the basis functions are defined as the solution of a local
elliptic partial differential equation, and are never explicitly
computed (ergo the name {\em virtual}) in the implementation of the method. 
Over each element $E$ in the mesh,
the trial and test functions belong to the local 
discretization (virtual) space that consist of polynomials
of order less than or equal to $k$ ($k$ is the order of the element) and in addition nonpolynomial functions. 
Since the virtual basis functions are unknown in each element, the VEM uses their elliptic polynomial projections to build the bilinear form (stiffness
matrix) and continuous linear functional (body force term) of the variational formulation.  Such projections are computable from the degrees of freedom within each element. The bilinear form on $E$ consists of two parts: the consistency term that approximates the stiffness matrix on a given polynomial space and the correction term that ensures stability. 
Essential boundary conditions in the VEM are imposed as in the
FEM, and element-level assembly procedures are used to form the
global stiffness matrix and force vector.

\smallskip
A notable advantage of the VEM is that computations can be done over meshes with arbitrarily-shaped convex and nonconvex 
(simple and nonsimple) elements without needing to compute the shape functions 
(generalized barycentric coordinates~\cite{Hormann:2017:GBC}) on such elements. In particular, {\em hanging nodes} on nonmatching (quadtree or weakly convex elements) meshes 
lead to conforming approximations. This facilitates modeling bimaterial interfaces, such as those that arise between polygonal inclusions and the matrix in cement-based composites, and also simplifies the imposition of contact conditions along interfaces~\cite{Wriggers:2016:VEM,Neto:2021:FPM}. Furthermore, mesh generation is simplified: an irregular-shaped stiff inclusion
can be modeled using a single {\em polygonal virtual element}~\cite{Artioli:2018:HOV}. Many of 
these and other positive attributes of the VEM have been emphasized and demonstrated in the virtual element literature, initially for low- and high-order
formulations for scalar elliptic problems~\cite{Ahmad:2013:EPV,Beirao:2014:HGV,Brezzi:2014:GBV,Beirao:2017:HOV,Dassi:2018:EHO}
and more recently for linear and nonlinear 
problems in the deformation of solid continua~\cite{Beirao:2013:VEL,Gain:2014:VEM,Wriggers:2016:VEM,Artioli:2017a:AO2,Chi:2017:SBF,Wriggers:2017:EVE,Mengolini:2019:EPV,Artioli:2020:CVE,Neto:2021:FPM,Sukumar:2021:VEA}.

In recent studies, the versatility of the VEM in composites modeling (multiphases, unit cell homogenization, and multiscale computations) has been shown~\cite{Pingaro:2019:FSH,Rivarola:2019:MAV,Cascio:2020:VEM,Rivarola:2020:VEI}.  In this paper, we provide 
comparisons of the VEM versus Voronoi-cell lattice models for modeling the elastic deformation of 
two-dimensional multiphase composites, which serve as a basis for modeling fracture in future work. The comparisons highlight attributes of the VEM, including its accuracy and flexibility in discretizing multiphase materials.

\section{Elastostatic Model for Multiphase Materials: Strong and Weak Formulations}\label{sec:bvp}
\begin{figure}
  \centering
   \includegraphics[scale=0.85]{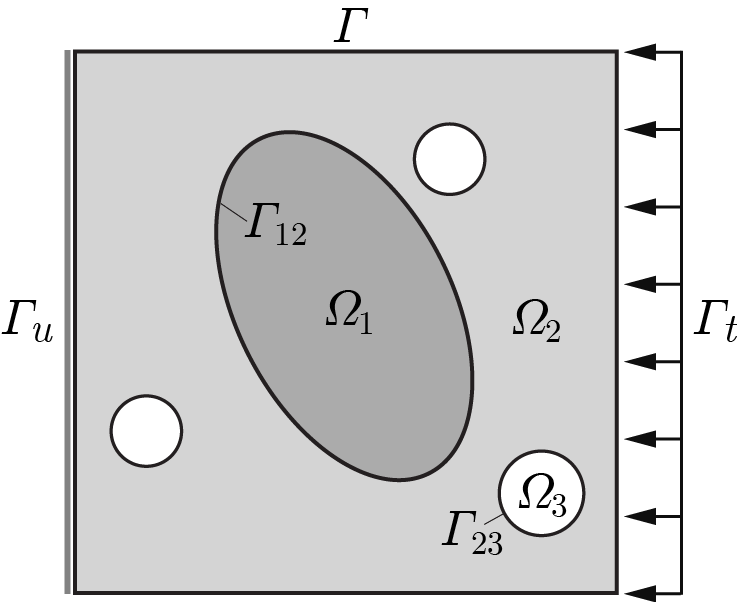}
  \caption{Elastostatic model problem for a three-phase composite.}
  \label{fig:3phase}
\end{figure}

Consider a 
linear elastic solid that occupies the domain $\Omega
\subset \Re^2$, with boundary $\Gamma = \partial \Omega$. The solid is composed
of $m$ isotropic, linearly elastic homogeneous materials, and the domain of each material is
$\Omega_i$, such that $\Omega = \overline{\Omega_1 \cup
\Omega_2 \cup \ldots \cup \Omega_m}$. The boundary that defines the material interface between $\Omega_i$
and $\Omega_j$ is denoted by $\Gamma_{ij}$. 
The material interface is assumed to be perfectly bonded.  
The external boundary $\Gamma = 
\overline{\Gamma_u \cup \Gamma_t}$, with
$\Gamma_u \cap \Gamma_t = \emptyset$. 
The boundary subsets $\Gamma_u$ and $\Gamma_t$ 
are where displacements and tractions are imposed, respectively. A
schematic of the model problem for a three-phase composite is shown in~\fref{fig:3phase}. 

In the absence of body forces, the governing equations of the
elastostatic boundary-value problem are:
\begin{subequations}\label{eq:bvp}
  \begin{align}
    \nabla \cdot \vm{\sigma} & = 0 \quad \textrm{in } \Omega,
     \label{eq:eqm} \\
    \vm{\sigma} &= \mathbb{C}_i : \vm{\varepsilon} \quad \textrm{in }\Omega_i,
      \label{eq:constitutive} \\
    \vm{\varepsilon} &= \nabla_s \vm{u} \quad \textrm{in } \Omega, 
      \label{eq:strain-disp}
  \end{align}
\end{subequations}
where $\vm{u}$ is the displacement field, $\nabla_s$ is the symmetric gradient
operator, $\vm{\varepsilon}$ is the small-strain tensor, $\vm{\sigma}$ is
the Cauchy stress tensor, and
$\mathbb{C}_i$ is the material moduli tensor for a homogeneous, linear elastic
isotropic material in the domain $\Omega_i$ ($i = 1, 2, \ldots, m$). The essential boundary
conditions, traction boundary conditions, and interface conditions are:
\begin{subequations}\label{eq:bc}
\begin{align}
    \vm{u} &= \bar{\vm{u}} \quad \textrm{on } \Gamma_u , \\
    \vm{n} \cdot \vm{\sigma} &= \bar{\vm{t}} \quad \textrm{on } \Gamma_t ,\\
    \left[\left[ \vm{u} \right]\right] &= \vm{0}
    \quad \textrm{on } \Gamma_{ij} , \label{eq:bimat-disp-equil} \\
    \left[\left[ \vm{n} \cdot \vm{\sigma} \right]\right] &= \vm{0}
      \quad \textrm{on } \Gamma_{ij} , \label{eq:bimat-traction-equil}
  \end{align}
\end{subequations}
where $\vm{n}$ is the unit vector that is normal to the indicated
boundary and $\left[\left[ \cdot \right]\right]$ is the jump operator that
represents the jump in its argument across the interface.

The presence of a material (weak) discontinuity is
met in the standard finite element method and also in the virtual element method via meshing the
domain and its internal and external boundaries so that $\Gamma_{ij}$ is
the union of element edges in the FE mesh. The additional advantage in
the VEM is that this condition is retained even if nodal conformity is not met on an edge that is shared by two elements ({\em hanging nodes} are allowed). However, unlike use of polygonal finite elements that require shape functions over weakly convex polygons (quadtree meshes)~\cite{Tabarraei:2008:EFE}, 
this is achieved in the VEM without the need to form the shape functions.

Let $\U{i}$ ($i = 1,2$) denote the affine subspace of functions in the Sobolev
space $H^1(\Omega)$ whose trace on $\Gamma_u$ is equal to $\bar{u}_i$ and whose
normal derivative on $\Gamma_{ij}$ is discontinuous. In addition, let
$\U{0}$ denote the linear subspace of functions in the Sobolev space
$H^1(\Omega)$ that vanish on $\Gamma_u$ and whose normal derivative on
$\Gamma_{ij}$ is discontinuous. The weak form
of~\eqref{eq:bvp} and~\eqref{eq:bc} is: find the trial displacement field
$\vm{u} \in \U{1} \times \U{2}$ such that
\begin{subequations}\label{eq:weak-form}
  \begin{align}
    \label{eq:weak-form-a}
    a ( \vm{u}, \vm{v} ) &=
      \ell ( \vm{v} ) \quad \forall \vm{v}
      \in \U{0} \times \U{0},  \\
    \intertext{where $\vm{v}$ is the test displacement field, and 
    the internal virtual work (bilinear form) $a(\cdot,\cdot)$ and the
    linear functional $\ell(\cdot)$ are given by}
    \label{eq:weak-form-b}
    a ( \vm{u}, \vm{v} ) := \int_\Omega & \vm{\sigma}(\vm{u}) :  \vm{\varepsilon}(\vm{v})
      \, d \vx  , \quad  \ell ( \vm{v} ) := \int_{\Gamma_t} \bar{\vm{t}} \cdot \vm{v} \, ds .
  \end{align}
\end{subequations}

\section{Virtual Element Method for Plane
         Elasticity}\label{sec:vem}
The formulation and implementation of the lowest-order
virtual element method for 2D and 3D 
solid continua is well-documented~\cite{Beirao:2014:HGV,Gain:2014:VEM,Artioli:2017a:AO2,Mengolini:2019:EPV,Sukumar:2021:VEA}. We follow the 
exposition in Sukumar and Tupek~\cite{Sukumar:2021:VEA} to present the main elements of VEM for 2D solid continua.

\subsection{Decomposition of the domain}\label{subsec:domain}
Let $\Omega \subset \Re^2$ be the problem domain and
${\cal T}^h$ a decomposition of $\Omega$ into nonoverlapping polygons 
(simple or nonsimple). The number of nodes in ${\cal T}^h$ is $N$. 
We refer to $E \in {\cal T}^h$
as an \emph{element}. 
The vertices of $E$ are denoted by $v_i$, and the coordinate of vertex 
$v_i$ by $\vx_i := (x_i,y_i)$. The diameter, centroid (barycenter) and area 
of $E$ are denoted by $h_E$. $\vx_E$ and $|E|$, respectively.
A polygon $E$ has $N_E$ vertices and $N_E$ edges, with the edges 
denoted by $e_i$ ($i=1,2,\ldots,N_E$). For the convergence proofs,  restrictions are placed on the shape-regularity of the
elements~\cite{Beirao:2013:BPV}.

\subsection{Polynomial spaces and virtual element space}
\label{subsec:VEspace}
Let $\mathbb{P}_k (E)$ be the function space on $E$ that consists of 
all polynomials of order less than or equal to $k$. By convention, $\mathbb{P}_{-1} = \{0\}$.  
The dimension of $\mathbb{P}_k(E)$ is denoted by $\textrm{dim}\, \mathbb{P}_k(E)$, and $\textrm{dim} \, \mathbb{P}_k(E) = (k+1)(k+2)/2$
in two dimensions.
The set consisting of the scaled monomials of order less than 
or equal to $k$ on $E$ is defined as ${\cal M}_k(E)$. In this paper, we 
use the first-order VEM ($k = 1$). In two dimensions,
\begin{equation*}
\mathbb{P}_1(E) = \left\{1, \ x , \ y \right\}, \ \ 
{\cal M}_1(E) = \left\{ 1 , \ \frac{x - x_E}{h_E} ,  \dfrac{y - y_E}{h_E} 
 \right\}
\end{equation*}
are the first-order (polynomial and scaled monomial, respectively)
basis sets. All elements in ${\cal M}_1(E)$ are of ${\cal O}(1)$.

\smallskip
Let $\vm{P}_1(E) = [\mathbb{P}_1(E)]^2$
be the 
polynomial basis for a vector field in $\Re^2$.
For planar linear elasticity, there are three
rigid-body (zero-energy) modes. Let
\begin{equation}\label{eq:xieta}
\xi := \dfrac{x - x_E}{h_E}, \quad \eta := \dfrac{y - y_E}{h_E} .
\end{equation}
For 2D solid continua, we define 
$\widehat{\vm{M}}(E) := \vm{M}_1(E)$ as the scaled monomial first-order vectorial
basis set:
\begin{equation}\label{eq:M}
\widehat{\vm{M}}(E) = 
\left[
\begin{Bmatrix} 1 \\ 0 \end{Bmatrix}, \ 
\begin{Bmatrix} 0 \\ 1 \end{Bmatrix}, \ 
\begin{Bmatrix} -\eta \\ \xi \end{Bmatrix}, \
\begin{Bmatrix} \eta \\ \xi \end{Bmatrix}, \
\begin{Bmatrix} \xi \\ 0 \end{Bmatrix}, \
\begin{Bmatrix} 0 \\ \eta \end{Bmatrix}
\right] , 
\end{equation}
where the first three vectorial bases in~\eqref{eq:M} contain 
the rigid-body modes. 

\smallskip
Let $V(E)$ denote the first-order virtual element space on element $E$.
The virtual element space for 2D solid continua is~\cite{Beirao:2013:VEL}:
\begin{align*}
V(E) = \Bigl\{  & \vm{v}^h : \vm{v}^h \in [H^1(E)]^2, \ 
        \Delta \vm{v}^h = \vm{0}, \\
         & \vm{v}^h |_e \in \vm{P}_1(e) \ \forall e \in \partial E, \
         \vm{v}^h |_{\partial E} \in [C^0]^2 (\partial E) 
          \Bigr\},
\end{align*}
where $\vm{v}^h$ is a piecewise continuous affine vector polynomial on the boundary of the polygon.

\subsection{Computation of energy projection matrices}
For linear elasticity,
we take the values of $\vm{v}^h$ at the vertices of the polygon as its degrees of freedom (DOFs).
Let $\{\phi_i\}_{i=1}^{N_E}$ be virtual canonical
basis functions that satisfy the Lagrange interpolation property,
$\phi_i(\vx_j) = \delta_{ij}$. 
Define the vectorial 
basis function matrix: 
\setcounter{MaxMatrixCols}{20}
\begin{equation} \label{eq:varphi}
\begin{split}
\vm{\varphi } 
&=
\begin{bmatrix}
\phi_1 & \dots &  \phi_{N_E}   &   0    &  \dots  &  0\\  
   0   & \dots &       0       &   \phi_1  &  \dots  &  \phi_{N_E}  
\end{bmatrix} \\
&:= 
\begin{bmatrix}
\vm{\varphi}_1  & \dots & \vm{\varphi}_{N_E} & 
\vm{\varphi}_{N_E+1} & \dots & \vm{\varphi}_{2N_E} 
\end{bmatrix},
\end{split}
\end{equation}
where $ \vm{\varphi}_{i} = \{ \phi_i \ 0 \}^T$ and
$\vm{\varphi}_{N_E+i} = \{ 0 \ \phi_i \}^T$ for $i = 1,\dots,N_E$, are the $2N_E$
vectorial basis functions. 
The trial displacement field in $E$ is:
\begin{equation}\label{eq:vh1} 
\vm{v}^h(\vx) = \sum_{i=1}^{2 N_E} \vm{\varphi}_i (\vx) v_i := 
\sum_{i=1}^{2 N_E} \vm{\varphi}_i (\vx) \texttt{dof}_i(\vm{v}^h),
\end{equation}
where $v_i$ are scalar coefficients and $\texttt{dof}_i(\cdot)$ extracts the $i$-th DOF of its argument.

Let $a_E^h(\cdot,\cdot)$ represent the discrete bilinear form of the countinuous operator in~\eqref{eq:weak-form-a}. The variational problem to determine the projector is determined via the energy orthogonality condition: 
\begin{subequations}\label{eq:orthog}
\begin{align}
    \label{eq:orthog-a}
    a_E^h(\bm{m}_\alpha,\bm{v}^h - \Pi^\varepsilon \bm{v}^h) &= 0 \ \ \forall \bm{m}_{\alpha} \in \bm{\widehat{M}}(E),
    \intertext{which is supplemented by the condition:}
    P_0(\bm{m}_\alpha,\bm{v}^h -\Pi^\varepsilon \bm{v}^h)  &= 0 \ \ (\alpha = 1,2,3), \\
    P_0(u,v) &= \frac{1}{N_E} \sum_{j=1}^{N_E} \vm{u}(\vx_j) \cdot \vm{v}(\vx_j),
\end{align}
\end{subequations}
where the projector $P_0(\cdot,\cdot)$ defines a discrete $L^2$ inner product on $E$. Note that
for $\alpha = 1,2,3$, \eqref{eq:orthog-a} yields $0 = 0$.

On using~\eqref{eq:varphi}, we define
\begin{equation}\label{eq:Pivarphi_i}
    \Pi^\varepsilon \vm{\varphi}_i = 
    \sum_{\beta = 1}^{6} \vm{m}_\beta  \pi_\beta^i 
    = \widehat{\vm{M}} \vm{\pi}^i \ \ 
    (i = 1,2,\dots,2N_E)
\end{equation}
as the projection of the $i$-th vectorial basis function onto the scaled monomial basis set, where
$\pi_\beta^i$ are unknown coefficients. 
On substituting $\vm{v}^h = \vm{\varphi}_i$ ($i = 1,2,\ldots 2N_E$)
in~\eqref{eq:orthog}, using the divergence theorem on the right-hand side and linear momentum balance
($\nabla \cdot \vm{\sigma} = \vm{0}$), we obtain the linear system of equations:
\begin{subequations}\label{eq:Pi}
\begin{align}
  \vm{G} \vm{\Pi} &= \tilde{\vm{B}},  \quad
  \vm{\Pi} = \vm{G}^{-1} \tilde{\vm{B}}, \\
  \vm{G}_{\alpha \beta} &= \begin{cases}
            \frac{1}{N_E} \sum_{j=1}^{N_E} \vm{m}_\alpha (\vx_j) \cdot \vm{m}_\beta (\vx_j) \ \ (\alpha=1,2,3) \\
            \vm{\sigma}(\vm{m}_\alpha ) : \vm{\varepsilon}(\vm{m}_\beta) \, |E| \ \ 
            (\alpha = 4, 5, 6)
           \end{cases} \! \! \! \! , \\ 
   \tilde{\vm{B}}_{\alpha i} &= \begin{cases}
              \frac{1}{N_E} \sum_{j=1}^{N_E} \vm{m}_\alpha (\vx_j) \cdot \vm{\varphi}_i (\vx_j) \ \  (\alpha=1,2,3) \\
              \label{eq:Pi-c}
              \vm{\sigma}(\vm{m}_\alpha) :
              \sum_{j=1}^{N_E} \int_{e_j} \! \vm{\varphi}_i \otimes \vm{n}_j \, ds \ \ (\alpha = 4, 5, 6)
           \end{cases} \! \! \! \! \! \! \! \! ,
\end{align}
\end{subequations}
where $\vm{\Pi} = [ \vm{\pi}^1, \vm{\pi}^2, \dots , \vm{\pi}^{2N_E} ]$ is the matrix representation of the projection of the canonical 
basis functions in
the scaled monomial basis set. The 
boundary integral in~\eqref{eq:Pi-c} can be exactly computed using a 
two-point Gauss-Lobatto quadrature scheme. 

Using the decomposition $\vm{\varphi}_i = \Pi^\varepsilon \vm{\varphi}_i + 
(1 - \Pi^\varepsilon ) \vm{\varphi}_i$ in the bilinear form and the orthogonality condition~\eqref{eq:orthog-a}, the stiffness
matrix can be expressed as:
\begin{align*}
\vm{K}_E &= \vm{K}_E^c + \vm{K}_E^s, \\
\vm{K}_E^c &= \vm{\Pi}^T  \tilde{\vm{G}} 
                            \vm{\Pi}, \\ 
    \vm{K}_E^s &= ( \vm{I} - \vm{D} \vm{\Pi} )^T \,
    a_E^h(\vm{\varphi},\vm{\varphi}) \
                     ( \vm{I} - \vm{D} \vm{\Pi} ), \\
                     \intertext{and to ensure stability we approximate $a_E^h(\vm{\varphi},\vm{\varphi})$ by a diagonal
                     matrix $\vm{S}_E^d$ that scales as $\vm{K}_E^c$:}
    \vm{K}_E^s &:= 
       ( \vm{I} - \vm{D} \vm{\Pi} ) \,  \vm{S}_E^d \,
       ( \vm{I} - \vm{D} \vm{\Pi} ),
\end{align*}
where $\vm{K}_E^c$ and $\vm{K}_E^s$ are the consistency and stabilization matrices, respectively, and
$\tilde{\vm{G}}$ is the matrix $\vm{G}$ with its first three rows set to zero. In addition,
$\vm{D}_{i\alpha} = \texttt{dof}_i (\vm{m}_\alpha)$ is the DOF-matrix and
the $i$-th diagonal entry of 
$\vm{S}_E^d$ is chosen as
$\max \bigl( \texttt{tr} \, (\vm{C})/3,  (K_E^c)_{ii} \bigr)$~\cite{Dassi:2018:EHO}, where $\vm{C}$ is the isotropic linear elastic constitutive matrix.

\subsection{Assembly and solution procedure}
On the natural boundary $\Gamma_t$, the virtual element shape functions are identical to piecewise linear finite elements. So on using~\eqref{eq:weak-form-b}, the
element force vector
$\vm{f}_E = \ell(\vm{\varphi})$ is computed. Having computed
$\vm{K}_E$ for each element, we then
perform standard finite element assembly procedures to form the global stiffness matrix $\vm{K}$ and the global force vector $\vm{f}$. On
incorporating the essential boundary conditions, the linear system is solved to obtain the nodal displacement vector $\vm{d}$.

\section{Voronoi-Cell Lattice Models}\label{sec:RBSN}

\begin{figure}[!tbh]
  \centering
\includegraphics[width=0.35\textwidth]{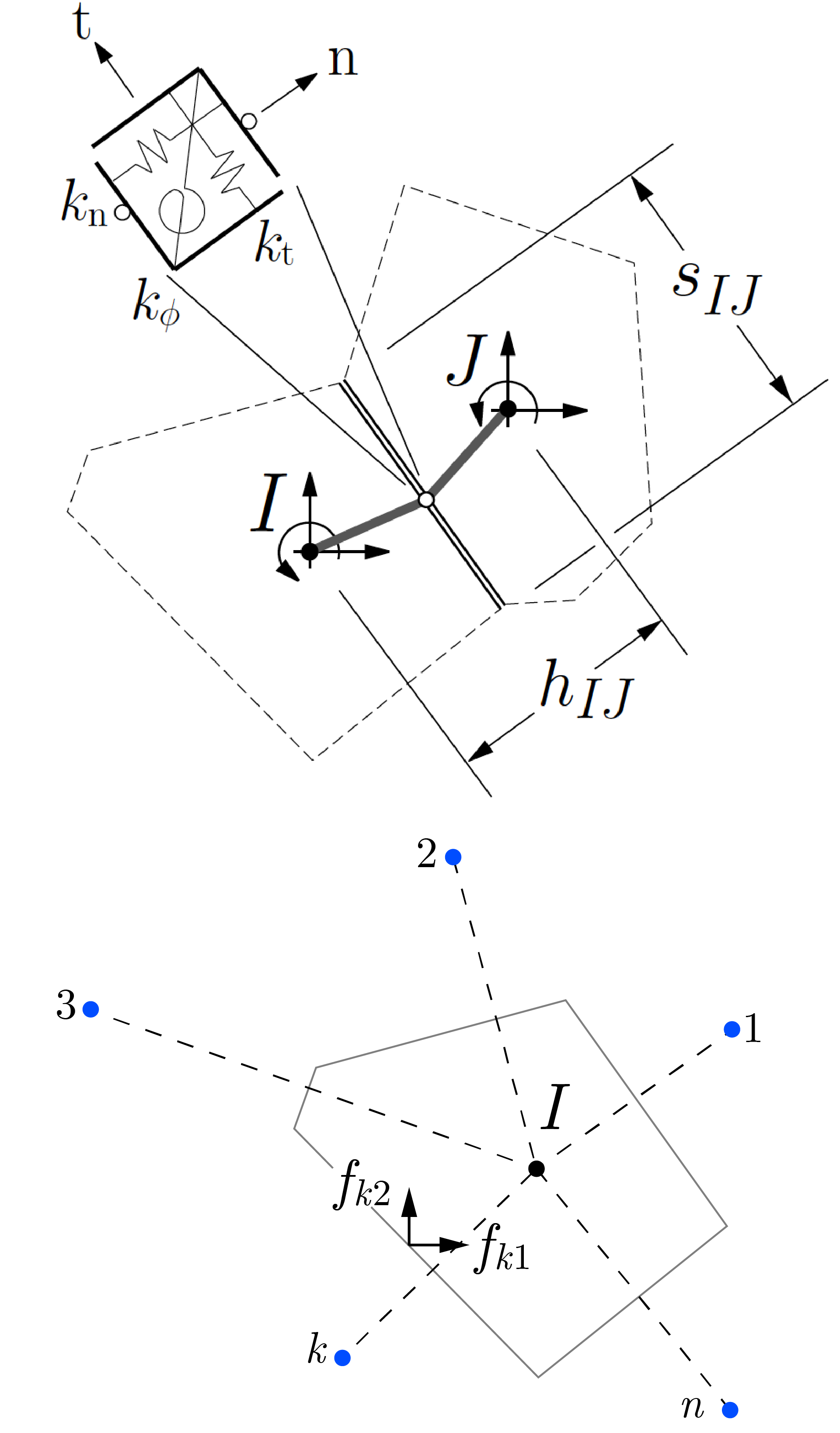}
  \caption{Voronoi-cell lattice element based on the rigid-body-spring concept (top);  and spring-set force components for nodal stress calculation (bottom).}
  \label{fig:rbsn}
\end{figure}

Particle-based lattice models share many features with classical lattice models, yet they differ in that each node is positioned within a geometric construct, or particle. In this case, the Voronoi diagram is used to partition the domain, such that each lattice node is associated with a Voronoi cell. For elasticity problems, the lattice element formulations are based on the rigid-body-spring concept of Kawai~\cite{Kawai:1978}. The Voronoi cells are assumed to be rigid and interconnected  via zero-size spring sets located midway along the facets common to neighboring cells (Fig.~\ref{fig:rbsn}). The stiffness matrices of these lattice elements are akin to those of ordinary frame elements and assemble into the system stiffness matrix in the conventional manner.
The spring sets contain components that are normal and tangential to the corresponding Voronoi facet. When these spring components have the same stiffness (i.e., $k_n = k_t$), the assembly of lattice elements is elastically homogeneous under uniform straining~\cite{Bolander:1998}, albeit with $\nu = 0$. 

The desired representation of both elastic constants (Young's modulus and Poisson's ratio) can be achieved, in a macroscopic sense, by appropriately setting the spring stiffness coefficients~\cite{Nagai:2005,Elias:2017,Elias:2020}. The assignment of spring coefficients depends on the type of loading. For the case of plane stress, and the type of random lattice considered herein, the coefficients can be determined according to~\cite{Elias:2020}
\begin{subequations}
\begin{align}
    \label{eq:pratio}
    \nu &= \frac{1-\alpha}{3+\alpha} , \\
    \label{eq:ymodulus}
    E &= E_0 \frac{2+2\alpha}{3+\alpha} ,
\end{align}
\end{subequations}
where $\alpha = k_t/k_n$ and $E_0$ is the effective elastic modulus at the element level.
Alternatively, Asahina 
et al.~\cite{Asahina:2015,Asahina:2017} have developed a procedure (based on $k_n = k_t$) in which the Poisson effect is introduced iteratively using the concept of auxiliary stress. This provides both local and global representations of elastic behavior.

\section{Numerical Examples}\label{sec:results}
Numerical simulations are performed using the VEM and VCLM models
on a suite of two-dimensional test problems.  Within each example, the same Voronoi tessellation is used to define each model.

\subsection{Patch test}\label{subsec:patch_test}
The VEM and VCLM models are assessed on the displacement patch test. Two forms of VCLM are considered in this example: (a) global representation of elastic behavior according to (\ref{eq:pratio}) and~(\ref{eq:ymodulus}); and (b) elastic behavior based on $k_n = k_t$ and the iterative introduction of the Poisson effect. 

\begin{figure*}[!htb]
  \centering
  \begin{subfigure}{0.27\textwidth}
      \includegraphics[width=\textwidth]{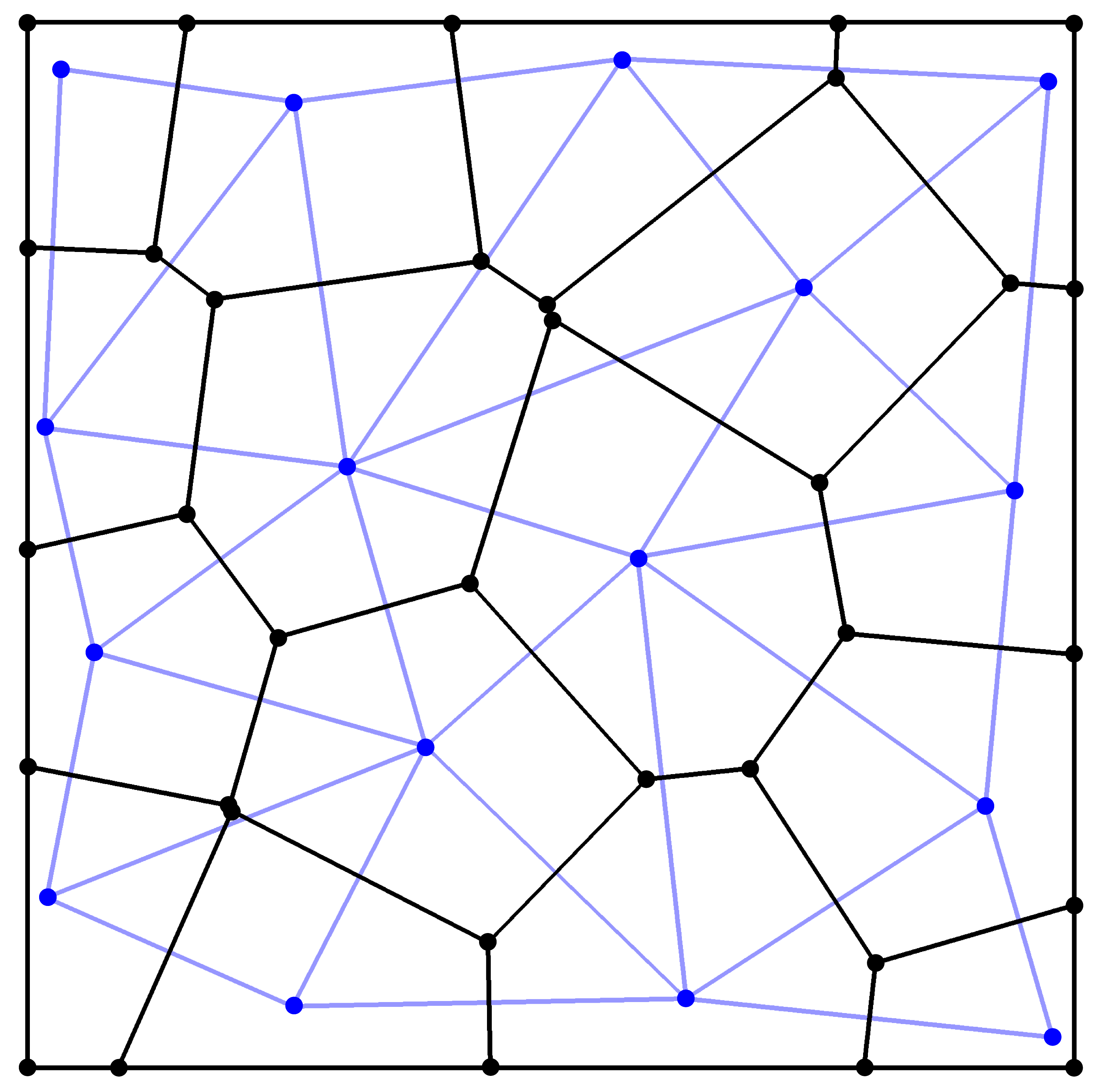}
      \caption{}\label{fig:meshes2D-a}
  \end{subfigure}
  \begin{subfigure}{0.27\textwidth}
      \includegraphics[width=\textwidth]{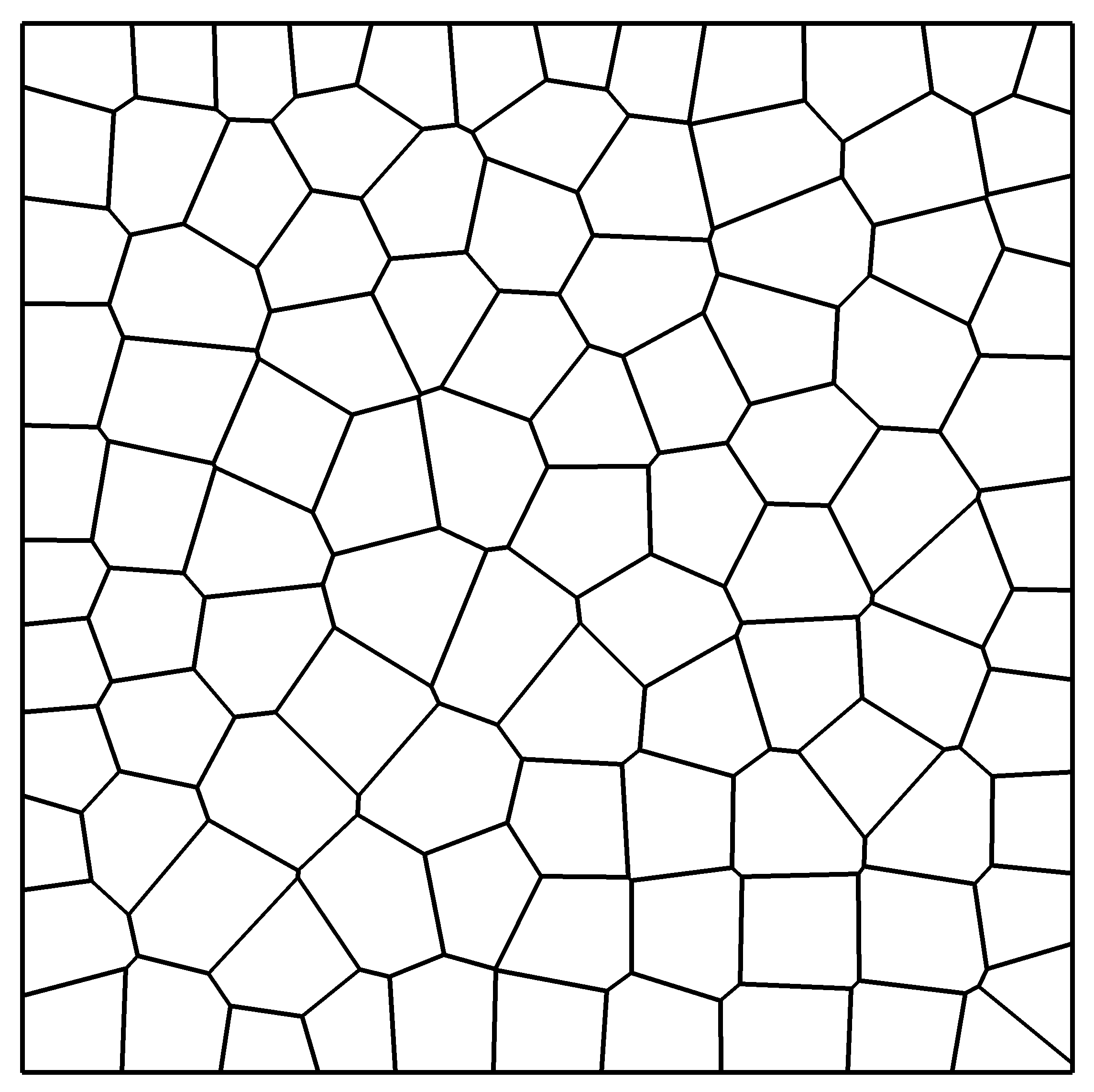}
      \caption{}\label{fig:meshes2D-b}
  \end{subfigure}
  \begin{subfigure}{0.26\textwidth}
      \includegraphics[width=\textwidth]{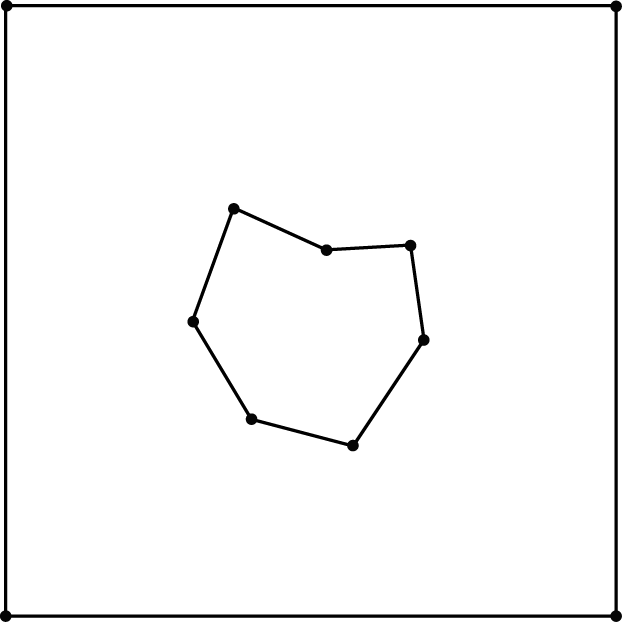}
      \caption{}\label{fig:meshes2D-c}
      \end{subfigure}
  \caption{Displacement patch test for different discretization schemes. The displacement field $u(x,y) = 1 + x + y$, $v(x,y) = 2 - 3x - 4y$ is imposed on nodes (or generator points) that lie on (or adjacent to) the boundary of the square. (a) Coarse Voronoi cell discretization, indicating element connectivities defined by the Voronoi vertices and lattice element connectivities defined by the Voronoi generator points; (b) fine Voronoi cell discretization; and (c) single virtual element discretization of a square region containing an irregularly-shaped inclusion. The mesh in (c) consists of two virtual elements: the inner element is a nonconvex septagon and the outer element is a nonsimply-connected polygon with an inner boundary (seven nodes in clockwise orientation) and an outer boundary (four nodes in counter-clockwise orientation).}\label{fig:meshes2D}
\end{figure*}

The meshes utilized for the patch test are shown in \fref{fig:meshes2D}. On all meshes, including
the mesh in~\fref{fig:meshes2D-c} that contains a nonsimply-connected element, the VEM passes the patch test as indicated in Table~\ref{t:L2norm}. The performance of the VCLM depends on the aforementioned assignment of its spring coefficients. When using~(\ref{eq:pratio}) and~(\ref{eq:ymodulus}), such that $k_n \neq k_t$, relative errors are ${\cal O}(10^{-2})$. 
This outcome has been viewed, with arguably some merit, as an effective means for representing the heterogeneity of concrete materials~\cite{Nagai:2005}. By setting $k_n=k_t$ and introducing the effect of Poisson's ratio using auxiliary stresses, however, the VCLM passes the patch test. For both VEM and VCLM ($k_n = k_t$), the errors in the stress components are also found to be within machine precision. We point out that if $k_n = k_t$ is used in VCLM, but without the auxiliary stress modifications in the
algorithm, then the displacement field matches the exact solution yet the  computed stresses correspond to the case of $\nu = 0$.

\begin{table}[htb]
\caption{$L^2$ norm of the error in the displacement for the
         patch test.}  \label{t:L2norm}
\begin{center}
\begin{tabular}{cccc} \hline  \\ [-2.5ex]
& \multicolumn{3}{c}{Discretization} \\ 
Model & Fig.\ \ref{fig:meshes2D-a} & Fig.\ \ref{fig:meshes2D-b} & 
  Fig.\ \ref{fig:meshes2D-c} \\[-2ex] 
& & & \\ \hline  
VEM & $2 \times 10^{-16}$ & $4 \times 10^{-16}$ & $1 \times 10^{-16}$ \\
VCLM ($k_n \neq k_t$)  & $4 \times 10^{-02}$ & $4 \times 10^{-02}$& --
\\ 
VCLM ($k_n = k_t$) & $3 \times 10^{-16}$ & $2 \times 10^{-16}$& --\\ 
\hline \\[-7mm]
\end{tabular}
\end{center}
\end{table}

\subsection{Bimaterial subjected to axisymmetric plane strain} \label{subsec:axisymm_inclusion}
\begin{figure}[!tbh]
  \centering
\includegraphics[width=0.48\textwidth]{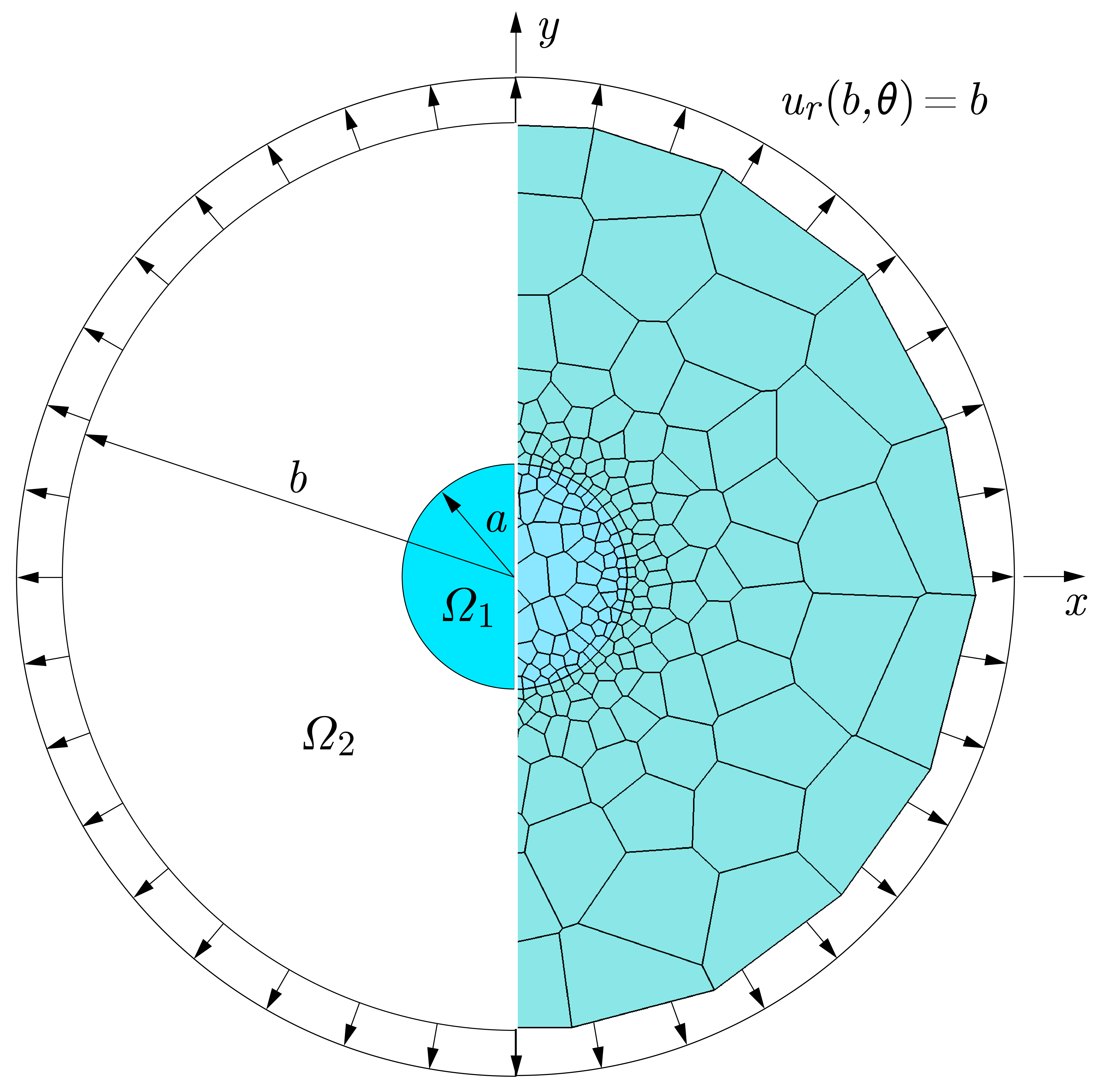}
\includegraphics[width=0.46\textwidth]{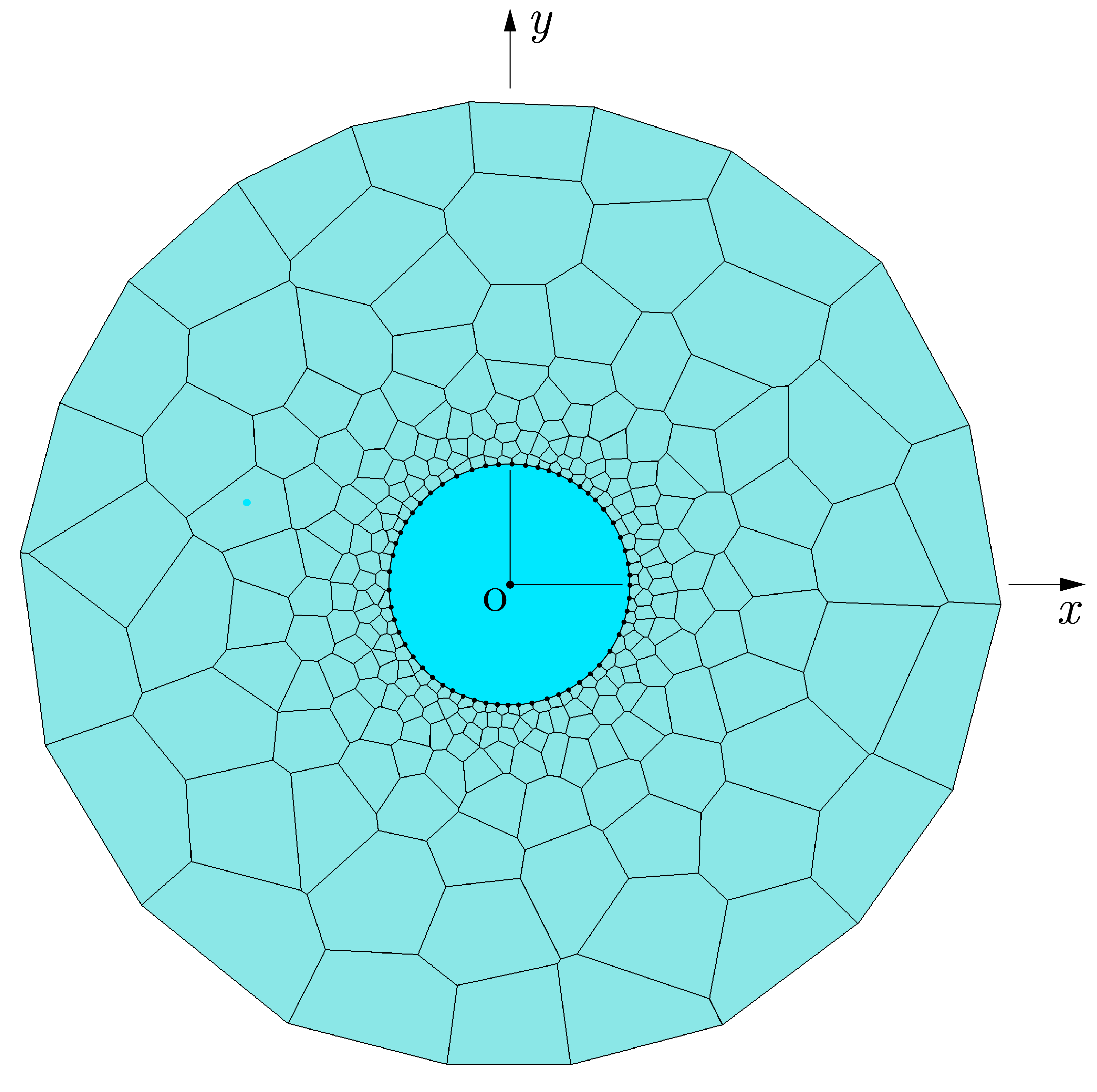}
  \caption{Bimaterial boundary-value problem. 
  Model problem and boundary conditions (top left),
  Voronoi-cell discretization of the domain (top right)
  and representation of central inclusion using a single element (bottom) are
  shown.}
  \label{fig:inclusion}
\end{figure}

We consider a two-phase composite that occupies a circular region of radius $b$. The inclusion, $\Omega_1$, is a disk of radius $a$ and the matrix is defined by
the region $\Omega_2 = \{r: a < r \le b \}$.  A radial displacement of
magnitude $b$ is imposed on $r = b$.
Due to axisymmetry, the displacement field is: $\vm{u}(r,\theta) = u_r(r) \vm{e}_r$ in which $u_r(r)$ is a nonzero radial displacement field.
A schematic illustration of the boundary-value problem is shown in~\fref{fig:inclusion}, along with a Voronoi-cell discretization of the bimaterial domain that defines both the VEM and lattice models.
This problem was first proposed in Sukumar et al.~\cite{Sukumar:2001:MHI}, and serves as a benchmark problem in computational solid mechanics~\cite{Schroder:2021:SBP}. 
The exact displacement field is~\cite{Sukumar:2001:MHI}:
\begin{align*}
    \begin{split}
    u_r (r) & = \begin{cases}
      \left[ \left( 1 - \frac{b^2}{a^2} \right) \alpha + \frac{b^2}{a^2} \right] r
        & 0 \leq r \leq a \\
      \left( r -\frac{b^2}{r} \right) \alpha + \frac{b^2}{r} &
        a < r \leq b
      \end{cases} \, , \\
    u_\theta & = 0 ,
  \end{split}
\end{align*}
where
\begin{equation*}
  \alpha = \frac{\left( \lambda_1 + \mu_1 + \mu_2 \right) b^2}
    {\left( \lambda_2 + \mu_2 \right) a^2 + \left( \lambda_1 + \mu_1 \right)
    \left( b^2 - a^2 \right) + \mu_2 b^2} ,
\end{equation*}
and $\lambda_1$ and $\mu_1$ and $\lambda_2$ and $\mu_2$ are the Lam\'{e}
parameters in $\Omega_1$ and $\Omega_2$, respectively.

This problem is solved for several values of the 
modular ratio $\eta$ = $E_1 / E_2$ with $a/b = 0.25$ and 
$\nu_1 = \nu_2 = 0.3$.
Table~\ref{t:L2normb} presents relative $L^2$ norm of the error in
the displacement field for the VEM and VCLM models.
For $\eta = 1$, the material is homogeneous, leading to uniform biaxial tension throughout the domain. Both methods simulate this condition with high precision.

\begin{table}[htb]
\caption{Relative $L^2$ norm of the error in the displacement for the bimaterial problem.}  \label{t:L2normb}
\begin{center}
\begin{tabular}{cccc} \hline  
{} & \multicolumn{3}{c}{Modular ratio $\eta$ } \\ 
{Model}      & 1 & 10 & 100\\ \hline \\[-2ex]
VEM & $1.2\times 10^{-15}$ & $1.5\times 10^{-3}$ & $1.8\times 10^{-3}$\\
VCLM & $2.5\times 10^{-10}$ & $2.7\times 10^{-3}$& $3.9\times 10^{-3}$\\
\hline
\end{tabular}
\end{center}
\end{table}

\begin{figure}[!tbh]
  \centering
    \includegraphics[width=0.48\textwidth]{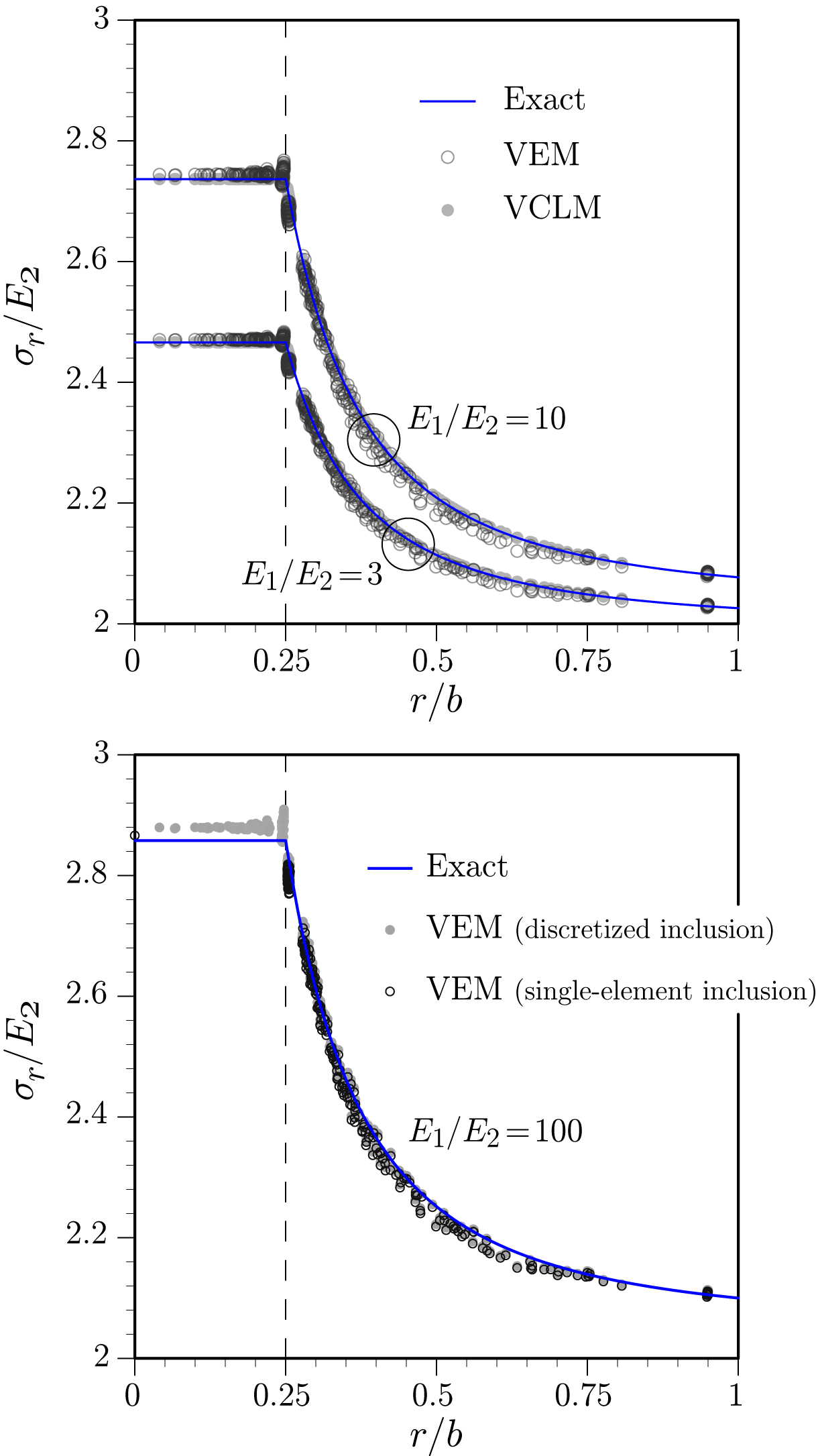}
  \caption{Radial profiles of radial stress for moderately stiff inclusions (top) and a stiff inclusion, highlighting the use of a single virtual element to represent the inclusion (bottom).}
  \label{fig:stressProfile}
\end{figure}

Figure~\ref{fig:stressProfile} plots major principal stress as a function of distance from the 
center of the inclusion. The 
exact solution for the stress components is:
\begin{subequations}
\begin{align}
    \sigma_{rr}(r) &= 2\mu_i \varepsilon_{rr}(r) + \lambda_i (\varepsilon_{rr}(r) + \varepsilon_{\theta \theta}(r) ) , \\
    \sigma_{\theta \theta}(r) &= 2\mu_i \varepsilon_{\theta \theta}(r) + \lambda_i (\varepsilon_{rr}(r) + \varepsilon_{\theta \theta}(r) ),
    \end{align}
\end{subequations}
where $\varepsilon_{rr}(r)$ = $du_r(r)/dr$ and $\varepsilon_{\theta \theta}(r)$ = $u_r(r)/r$; the subscript on the Lam\'{e} constants indicates the material subdomain.

For the lattice simulations, the volume-averaged stresses are computed at the lattice 
nodes via the relation~\cite{Bardet:2001}
\begin{equation}\label{eq:average_stress}
\begin{split}
   \bar{ \bm{\sigma}}  &= \frac{1}{V} \int_C \vm{\sigma} \, d\vx
   = \frac{1}{V} \int_C \bigl[ \vm{I} \cdot \vm{\sigma} + 
      \vx \otimes (\nabla \cdot \vm{\sigma} ) \bigr] \, d\vx \\
   &= \frac{1}{V} \int_C ( \vx \otimes \vm{\sigma}^T ) \cdot \nabla \, d\vx = \frac{1}{V} \int_{\partial C} \vx \otimes ( \vm{n} \cdot \vm{\sigma} )
     \, dS \\
   &= \frac{1}{V} \int_{\partial C} \vx \otimes \vm{t} \, dS 
   \approx \frac{1}{V} \sum_{k=1}^{n_f} 
       \vx_k \otimes \vm{f}_k ,
\end{split}
\end{equation}
where $C$ is the Voronoi cell and $\partial C$ is its boundary. 
The stress tensor is assumed to be divergence-free (no body forces are present) and $\vm{t}$ are the boundary
tractions. In addition, 
$\vm{f}_k$ is a system of $n_f$ external forces acting on the corresponding Voronoi cell having volume $V$; the forces act at locations $\vm{x}_k$ with respect to the cell node. Note that the stress tensor as defined in~\eqref{eq:average_stress} is not symmetric, which is consistent
with the behavior of a discrete lattice
model as a micropolar (Cosserat) continuum~\cite{Bardet:2001}.

The simulated stress profiles in the radial direction (Fig.~\ref{fig:stressProfile}, top) agree well with theory for the range of $\eta$ values considered. For the VCLM model, the elements that span the $\Gamma_{12}$ boundary are assigned the harmonic mean values of the properties of domains $\Omega_1$ and~$\Omega_2$. Furthermore, the Voronoi generator points that define the boundary are positioned close to the boundary. These conditions improve the accuracy of the stress calculations near the boundary.

For cases of dispersed stiff inclusions, where the degree of modular mismatch is high, the stress field in the inclusions is approximately uniform. Such inclusions can be represented using a single VEM element, as shown in Fig.~\ref{fig:inclusion} (bottom). The radial stress value calculated for this single element, and plotted at $r/b$~=~0 in Fig.~\ref{fig:stressProfile} (bottom), has a relative error of $3.0 \times 10^{-3}$. For the region $r > a$, the radial stress values for the two cases (i.e., for the fully discretized and single-element representations of the inclusion) are essentially the same.

\subsection{Three-phase composite} \label{subsec:3phase}
\begin{figure*}[!tbh]
  \centering
\includegraphics[width=0.99\textwidth]{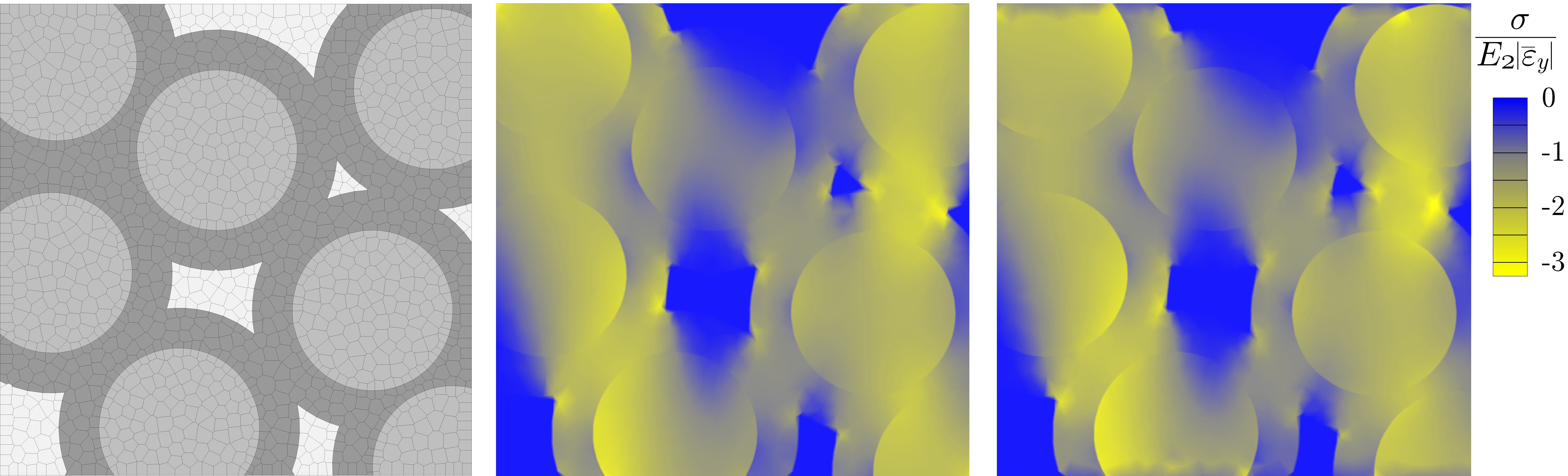}
  \caption{Model of macro-porous concrete under vertical compressive loading. 
  Domain discretization (left), and contours of minor
  principal stress using VEM (middle) and VCLM (right) are shown.}
  \label{fig:poc}
\end{figure*}

Capabilities of the VEM, and its correspondence to the VCLM, are further demonstrated through elastic analysis of a three-phase composite material. Figure~\ref{fig:poc}a shows the planar discretization of a model porous concrete, in which disk-shaped aggregate inclusions are coated with a uniformly thick layer of hardened cement paste; the lightest shaded regions represent large-scale porosity between the paste layers. Compressive load is applied in the form of a uniform downward displacement of the uppermost vertices (or nodes), producing an average vertical strain of $\overline{\epsilon}_y$. The modular ratio of the inclusion and cement paste materials is $E_1/E_2$ = 3. Poisson's ratio is set as $\nu$ = 0.2 for both phases, and plane stress conditions are assumed.

For the applied loading, Fig.~\ref{fig:poc} shows contours of minor principal stress, which highlight the nonuniform transfer of load through the material. Stress risers occur due to the stiff inclusions and large-scale porosity. The corresponding results for the VCLM are quite similar, except for differences that appear along the loaded faces of the models. These differences arise from the assignment of the displacement boundary conditions (at either the Voronoi vertices or generator points for the VEM and VCLM, respectively) and the difficulties in calculating nodal stress in the VCLM along the constrained boundaries.

\section{Conclusions}\label{sec:conclusions}
The challenges in simulating the mechanical behavior of multiphase composite materials include the effective, accurate modeling of elastic behavior, which is a determining factor for nonlinear material behavior. In this paper, we have
investigated the use of the virtual element method (VEM) for modeling the deformation of such composite materials. Displacements and element stress values were compared with theory and those of Voronoi-cell lattice models (VCLM) based on the same (dual) discretization scheme. In practical terms, both the VEM and VCLM approaches provided 
comparably 
accurate results. However, the VCLM required an iterative procedure to satisfy the patch test (elastic homogeneity) for arbitrary Poisson's ratio. In addition, VCLM nodes reside within the material domain, rather than on the domain boundaries, which complicates domain discretization and the assignment of boundary conditions. In this sense, the VEM has significant advantages. Furthermore, the VEM allows for stiff inclusions to each be modeled using a single polygonal element, which simplifies meshing relative to other approaches including the finite element method, particularly for irregularly shaped inclusions. A promising direction of future work is to use recent advances in the VEM on mesh-independent modeling of cracks~\cite{Benvenuti:2021:EVE} to simulate the deformation of multiphase composite materials, including cement-based composites, and their transition from continuous to discontinuous behavior.

\end{document}